\documentclass{amsart}

\newtheorem{theorem}{Theorem}[section]

\theoremstyle{definition}

\theoremstyle{remark}

\numberwithin{equation}{section}

\begin{document}

\title{The Mathematical Universe in a Nutshell}

\author{K. K. Nambiar}
\address{Formerly, Jawaharlal Nehru University, New Delhi, 110067, India}
\curraddr{1812 Rockybranch Pass, Marietta, Georgia, 30066-8015}
\email{kannan@rci.rutgers.edu}


\date{July 12, 2002}

\commby{}


\begin{abstract}
The mathematical universe discussed here gives models 
of possible structures our physical universe can have.
\end{abstract}

\maketitle



\section{Introduction}

\noindent 
When we talk about the universe, we usually mean the physical universe 
around us, but then we must recognize that we are capable of visualizing 
universes which are quite different from the one we live in. The 
possible universes that we can \emph{reasonably} imagine is what we 
collectively call the \emph{mathematical universe}. 
There is yet another universe which we may call the \emph{spiritual 
universe}, a universe hard to avoid and to define. We will accept a 
weak definiton of the spiritual universe as the collection of 
unverifiable persistent \emph{emotional beliefs} within us 
that are difficult to analyze. 
Our description makes it clear that spiritual universe is 
beyond and the physical universe is within the mathematical universe. 
As a matter of record, here are the three universes we are 
interested in:
\begin{itemize}
\item
Physical universe 
\item
Mathematical universe
\item
Spiritual universe
\end{itemize}
The purpose of this paper is to discuss the mathematical 
universe in some detail, so that we may have a deeper 
understanding of the physical universe and a greater appreciation 
for the spiritual universe.

\section{Mathematical Logic}

We want to talk about what the mathematicians call an axiomatic 
derivation and what the computer scientists call a string manipulation.
A familiar, but unemphasized fact is that it is the string 
manipulation of the four Maxwell's equations that allowed us 
to predict the radiation of electromagnetic waves from a dipole 
antenna. Similarly, it is the derivations in Einstein's theory 
of relativity that convinced us about the bending of light rays 
due to gravity. It is as though nature dances faithfully to the 
tune that we write on paper, as long as the composition strictly 
conforms to certain strict mathematical rules. 

It was the Greeks who first realized the power of the axiomatic 
method, when they started teaching elementary geometry to their 
school children. Over the years logicians have perfected the method 
and today it is clear to us that it not necessary to draw  
geometrical figures to prove theorems in geometry. Without belaboring 
the point, we will make a long story short, and state some of the 
facts of mathematical logic that we have learned to accept.

Mathematical logic reaches its pinnacle when it deals with 
Zermelo-Fraenkel set theory (ZF theory). Here is a set of 
axioms which define ZF theory.
\begin{itemize}
\item
Axiom of Extensionality: $\forall x(x\in a \Leftrightarrow x\in b) 
\Rightarrow (a=b)$. \\
Two sets with the same members are the same.
\item
Axiom of Pairing: $\exists x\forall y(y\in x \Leftrightarrow 
y=a\vee y=b)$. \\
For any $a$ and $b$, there is a set $\{a,b\}$.
\item
Axiom of Union: $\exists x \forall y (y \in x \Leftrightarrow 
\exists z \in x(y \in z))$. \\
For any set of sets, there is a set which has exactly the 
elements of the sets of the given set. 
\item
Axiom of Powerset: $\forall x \exists y (y \in x \Leftrightarrow 
\forall z \in y (z \in a))$. \\
For any set, there is a set which has exactly the subsets of the 
given set as its elements.
\item 
Axiom of Infinity: $\exists x (\emptyset \in x \wedge \forall 
y \in x (y \cup \{y\} \in x))$. \\
There is a set which has exactly the natural numbers as its elements.
\item 
Axiom of Separation: $\exists x \forall y 
(y \in x \Leftrightarrow y \in a \wedge A(y))$. \\
For any set $a$ and a formula $A(y)$, there is a set containing 
all elements of $a$ satisfying $A(y)$.
\item 
Axiom of Replacement: $\exists x \forall y \in a 
(\exists z A(y,z) \Rightarrow \exists z \in x A(y,z))$. \\ 
A new set is created when every element in a given set is replaced 
by a new element.
\item 
Axiom of Regularity: $\forall x (x \neq \emptyset \Rightarrow 
\exists y (y \in x \wedge x \cap y = \emptyset)) $.\\
Every nonempty set $a$ contains an element $b$ such that 
$a \cap b = \emptyset$.
\end{itemize}
These axioms assume importance when they are applied to infinite 
sets. For finite sets, they are more or less obvious. The real 
significance of these axioms is that they are the only strings, 
other than those of mathematical logic itself, 
that can be used in the course of a proof in 
set theory. Here is an example of a derivation using the axiom of 
regularity, which saved set theory from disaster.
\begin{theorem}
$a \notin a$.
\end{theorem}
\noindent 
Proof: $a \in a$ leads to a contradiction as shown below.
\begin{quote}
$a\in a \Rightarrow a \in [\{a\}\cap a]$ ... (1).\\ 
$b \in \{a\} \Rightarrow b = a$ ... (2). \\
Using axiom of regularity and (2),\\
$\{a\} \cap a = \emptyset$, which contradicts, 
$a \in [\{a\}\cap a]$ in (1).
\end{quote}
Even though we are in no position to prove it, 
over a period of time we have built up enough confidence in 
set theory to believe that there are no contradictions in it. 

\section{Generalized Anthropic Principle and The Book}

It is generally accepted that all of mathematics can be 
described in terms of the concepts of set theory, which in 
turn means that we can, in principle, axiomatize any 
branch of science, if we so wish. This allows us to 
conceive of an axiomatic theory which has all of known
science in it, and all the phenomena we observe in the 
universe having corresponding derivations. 
Since every derivation in a theory can be considered 
as a well-formed formula, we can claim that the set of derivations 
in our all-encompassing theory can be listed in the lexicographic 
order, with formulas of increasing length. A book which lists 
\emph{all} the proofs of this all-encompassing mathematics is 
called \emph{The Book}. The concept of The Book is an invention 
of the mathematician, Paul Erd\"os, perhaps the most prolific 
mathematician of the twentieth century. Note that the book 
contains an infinite number of proofs, and also that the 
lexical order is with respect to the proofs and not with 
respect to theorems. It was this book that David Hilbert, the 
originator of formalism, once wanted to rewrite with theorems 
in the lexical order, which of course, turned out to be an 
unworkable idea. 

Note that a computer can be set up
to start writing the book. We cannot expect the computer to stop, since 
there are an infinite number of proofs in our theory. Thus, a computer 
generated book will always have to be unfinished, the big difference 
between a computer generated book and The Book is that it is a 
\emph{finished} book. 

When discussing cosmology, a notion that is often invoked is 
called the \emph{anthropic principle}.
The principle states that we see the universe the way it is, 
because if it were any different, we would not be here to see 
it \cite{Hawking:UN}.
We generalize this concept as follows.
\begin{quote} 
\emph{Generalized Anthropic Principle\/:} Every phenomenon in 
the universe has a corresponding derivation in The Book.
\end{quote}
Note that the generalized anthropic principle does not claim that 
there is a phenomenon corresponding to every derivation. Such 
derivations are part of mathematics, but not of physics, in other 
words, we consider the physical universe as part of the mathematical 
universe.

\section{Expanding Universe}

As a preliminary to the understanding of the expanding universe,
we will first talk about a perfectly spherical balloon whose radius 
is increasing with velocity $U$, starting with $0$ radius and 
an elapsed time $T$. If we write $UT$ as $R$, we have the 
volume of the balloon as $(4/3)\pi R^{3}$, 
surface area as $4\pi R^{2}$ and the length of a 
great circle as $2\pi R$. In our expanding balloon
it is easy to see that two points which are a distance $R\theta$ 
apart from each other will be moving away from each other with 
velocity $U\theta$. Since $(R/U)=T$, it follows that a measurement 
of the relative movement of two spots on the balloon will allow us
to calculate the age of the balloon.

The facts about the expanding universe is more or less like that 
of the balloon, except that instead of a spherical surface, we 
have to deal with the hyper surface of a $4$-dimensional sphere
of radius $R$. If we use the same notations as before, we have the 
hyper volume of the hyper sphere as $(\pi^{2}/2)R^{4}$, the hyper 
surface as $2\pi^{2}R^{3}$, and the length of a great circle as 
$2\pi R$. We can calculate the age of the universe by measuring the
velocity of a receding galaxy near to us. If we assume the velocity 
of expansion of the hyper sphere as $U$, we have $UT=R$ and 
the volume of the universe as $2\pi^{2}R^{3}$.

We would have been more realistic, if we were to start off our 
analysis with a warped balloon, but then our intention here 
is only to discuss what is mathematically possible and not what 
the actual reality is.

\section{Intuitive Set Theory}

ZF theory which forms the foundations of mathematics gets 
simplified further to Intuitive Set Theory (IST),
if we add two more axioms to it as given below \cite{Nam:VIST}. 
If $k$ is an ordinal, we will write ${{\aleph_{\alpha}}\choose{k}}$ 
for the cardinality of the set of all subsets of
${\aleph_{\alpha}}$ with the same cardinality as $k$.

\textbf{Axiom of Combinatorial Sets:}
\[ 
\aleph_{\alpha+1}= 
{{\aleph_{\alpha}}\choose{\aleph_{\alpha}}}. 
\]

We will accept the fact that every number in the interval $(0,1]$ can be 
represented \textit{uniquely} by an 
\textit{infinite} nonterminating binary sequence. 
For example, the infinite binary sequence 
\[
.10111111 \cdots
\]
can be recognized as the representation for the 
number $3/4$ and similarly for other numbers.
This in turn implies that an \emph{infinite recursive} subset of positive 
integers can be used to represent numbers in the interval $(0,1]$.
It is known that the cardinality of the set $R$ of such recursive subsets 
is $\aleph_{0}$. Thus, every $r\in R$ represents a real number 
in the interval $(0,1]$.

We will write 
\[
{{\aleph_{\alpha}}\choose{\aleph_{\alpha}}}_{\!\!r},
\]
to represent the cardinality of the set of all those subsets of
${\aleph_{\alpha}}$ of cardinality ${\aleph_{\alpha}}$ 
which contain $r$, and also write 
\[
\Biggl \{{{\aleph_{\alpha}}\choose{\aleph_{\alpha}}}_{\!\!r}~
\bigg |~r\in R~ \Biggr \}
= {{\aleph_{\alpha}}\choose{\aleph_{\alpha}}}_{\!\!R}.
\]
We will define a \emph{bonded sack} as a collection which can appear 
only on the left side of the binary relation $\in$ and not on the 
right side. What this means is that a bonded sack has 
to be considered as an integral unit from which not even the 
axiom of choice can pick out an element. For this reason, we 
may call the elements of a bonded sack \emph{figments}.

\textbf{Axiom of Infinitesimals:}
\[
(0,1]={{\aleph_{\alpha}}\choose{\aleph_{\alpha}}}_{\!\!R}.
\]
The axiom of infinitesimals makes it easy to visualize the 
unit interval $(0,1]$.

We derive the generalized continuum hypothesis from the 
axiom of combinatorial sets as below: 
\[
2^{\aleph_{\alpha}}= 
{{\aleph_{\alpha}}\choose{0}}+ 
{{\aleph_{\alpha}}\choose{1}}+
{{\aleph_{\alpha}}\choose{2}}+
\cdots 
{{\aleph_{\alpha}}\choose{\aleph_{0}}}+
\cdots
{{\aleph_{\alpha}}\choose{\aleph_{\alpha}}}.
\]
Note that ${{\aleph_{\alpha}}\choose{1}}={{\aleph_{\alpha}}}$. 
Since, there are $\aleph_{\alpha}$ terms in this addition and 
${{\aleph_{\alpha}}\choose{k}}$ 
is a monotonically nondecreasing function of $k$, 
we can conclude that 
\[
2^{\aleph_{\alpha}}= 
{{\aleph_{\alpha}}\choose{\aleph_{\alpha}}}.
\]
Using axiom of combinatorial sets, we get 
\[
2^{\aleph_{\alpha}}= 
\aleph_{\alpha+1}.
\]

The concept of a bonded sack is significant 
in that it puts a limit beyond which the interval $(0,1]$ cannot 
be pried any further.
The axiom of infinitesimals allows us to visualize the unit 
interval $(0,1]$ as a set of bonded sacks, with cardinality $\aleph_{0}$. 
Thus, $\tbinom{\aleph_{\alpha}}{\aleph_{\alpha}}_{\!r}$ 
represents an \emph{infinitesimal} or \emph{white hole} 
or \emph{white strip} corresponding to the number $r$ in the 
interval $(0,1]$.

\section{Universal Number System}

A real number in the binary number system is usually defined as 
a two way binary sequence around a binary point, written as 
\[
xxx.xxxxx\ldots
\]
in which the left sequence is finite and the right sequence is 
nonterminating. Our discussion 
earlier, makes it clear that the concept of a real number and a 
white strip are equivalent.
The two way infinite sequence we get when we flip the real 
number around the binary point, written as 
\[
\ldots xxxxx.xxx
\]
we will call a \emph {supernatural number} or a \emph {black stretch}.
The set of white strips we will call the \emph{real line} and the 
set of black stretches the \emph{black whole}. 
Since there is a one-to-one correspondence between the white strips 
and the black stretches,
it follows that there is a duality between the real line and the 
black whole. 

The name black stretch is supposed to suggest that it 
can be visualized as a set of points distributed 
over an infinite line, but it should be recognized
as a bonded sack, which the axiom of choice cannot access. Our 
description of the black whole clearly indicates that it can be used 
to visualize what is beyond the finite physical space.

\section{Conclusion}

We will conclude with a few remarks about mathematical logic, 
which give us some indication why we cannot afford to ignore 
the spiritual universe. 
G\"odel tells us that there is no logical way to 
establish that there are no contradictions in ZF theory, which 
forms the foundations of mathematics. We are confident about 
our mathematics only because it has worked well for us for the 
last two thousand years. 
Since, any set of axioms is a set of beliefs, it follows 
that any theory is only a set of beliefs. Since, any individual is  
the sum total of h(is)er beliefs (axioms) and rational thoughts 
(derivations), no individual, including scientists, can claim 
to be infallible on any subject matter. If an honest scientist is 
called to appear in the ultimate court of nature, (s)he can use 
The Book for taking the oath, and the most (s)he can say 
is: \emph{I solemnly swear that if I am sane, I will tell nothing 
but the truth, but never the whole truth}.

\bibliographystyle{amsplain}

\end{document}